\documentclass[reqno]{amsart}
\usepackage{amssymb,amsmath}

\theoremstyle{plain}
\newtheorem{theorem}{Theorem}

\newtheorem{lemma}{Lemma}

\theoremstyle{definition}
\newtheorem{definition}{Definition}
\theoremstyle{remark}

\numberwithin{equation}{section}
\newcommand{\B}{\boldsymbol}
\newcommand{\Bf}[1]{\ensuremath{\B{#1}}}
\newcommand{\Bfd}[2]{\ensuremath{\B{#1}_{#2}}}
\newcommand{\Bfu}[2]{\ensuremath{\B{#1}^{#2}}}
\newcommand{\Bfdu}[3]{\ensuremath{\B{#1}_{#2}^{#3}}}
\newcommand{\Lfd}[2]{\ensuremath{{#1}_{#2}}}
\newcommand{\Lf}[1]{\ensuremath{#1}}
\newcommand{\Lfu}[2]{\ensuremath{{#1}^{#2}}}
\newcommand{\Lfdu}[3]{\ensuremath{{#1}_{#2}^{#3}}}

\newcommand{\Fdim}[1]{\mathbb{F}_{#1}}
\newcommand{\Rdim}[1]{\mathbb{R}^{#1}}
\newcommand{\Rdimp}[1]{\mathbb{R}^{#1}_{+}}
\newcommand{\Sdim}[1]{\mathbb{S}_{#1}}

\newcommand{\Sdimp}[1]{\mathbb{S}_{#1}^{+}}

\newcommand{\ScrOdim}[1]{\mathcal{O}({#1})}
\newcommand{\ScrUdim}[1]{\mathcal{U}({#1})}
\newcommand{\seq}[2]{\ensuremath{{#1}_{1},\ldots,{#1}_{#2}}}

\newcommand{\Sumil}[1]{\sum_{i=1}^{#1}}

\newcommand{\alf}{\ensuremath{\alpha}}
\newcommand{\bet}{\ensuremath{\beta}}
\newcommand{\eps}{\ensuremath{\epsilon}}

\newcommand{\gam}{\ensuremath{\gamma}}

\newcommand{\lam}{\ensuremath{\lambda}}
\newcommand{\Lam}{\ensuremath{\Lambda}}

\newcommand{\om}{\ensuremath{\omega}}

\newcommand{\sig}{\ensuremath{\sigma}}

\def\nmsp{\negmedspace}
\def\La{\mathcal{L}}

\def\Fkk{\Fdim{k\times k}}
\def\Fnk{\Fdim{n\times k}}
\def\Fkk{\Fdim{kk}}
\def\Fnk{\Fdim{nk}}

\def\FnkC{\Fnk^{\mathbb{C}}}
\def\Ok{\ScrOdim{k}}
\def\On{\ScrOdim{n}}
\def\Uk{\ScrUdim{k}}
\def\Un{\ScrUdim{n}}

\def\Rlp{\Rdimp{1}}
\def\Rn{\Rdim{n}}
\def\Rk{\Rdim{k}}
\def\Rkp{\Rdimp{k}}
\def\Sk{\Sdim{k}}

\def\Skp{\Sdimp{k}}

\def\Sumilk{\Sumil{k}}
\def\Sumiln{\Sumil{n}}

\def\ca{\Lfd{c}{1}}
\def\ck{\Lfd{c}{k}}
\def\csq{\Lfu{c}{2}}
\def\cst{\Lfu{c}{*}}
\def\cstsq{\Lfu{c}{*2}}
\def\casq{\Lfdu{c}{1}{2}}
\def\cksq{\Lfdu{c}{k}{2}}
\def\cg{\Lfd{c}{g}}
\def\cphibu{\Lfd{c}{\phi}(\cdot)}
\def\cphi{\Lfd{c}{\phi}}

\def\eilsq{\Lfdu{e}{(i,\lam)}{2}}
\def\gbu{g(\cdot)}

\def\ksq{\Lfu{k}{2}}
\def\vjj{\Lfd{v}{jj}}
\def\waa{\Lfd{w}{11}}

\def\wii{\Lfd{w}{ii}}
\def\wij{\Lfd{w}{ij}}
\def\wjj{\Lfd{w}{jj}}
\def\wkk{\Lfd{w}{kk}}
\def\Bc{\Lfd{B}{c}}
\def\Bcc{\Lfdu{B}{c}{c}}
\def\Clam{\Lfd{C}{\lam}}
\def\DZ{\Lfd{D}{\!Z}}
\def\DFlam{\Lfd{D\!F}{\!\lam}}
\def\Fa{\Lfd{F}{\!1}}
\def\Fc{\Lfd{F}{\!c}}
\def\Go{\Lfd{G}{\!0}}
\def\GCVlam{G\!C\!V_{\!\lam}}
\def\HKBlam{H\!K\!B_{\!\lam}}
\def\PRElam{P\!R\!E\!S\!S_{\!\lam}}
\def\Rlam{R(\!\lam)}

\def\Rsq{\Lfu{R}{2}}
\def\Sc{\Lfd{S}{c}}
\def\Uisq{\Lfdu{U}{i}{2}}
\def\Va{\Lfd{V}{1}}
\def\Vb{\Lfd{V}{2}}
\def\Vk{\Lfd{V}{k}}
\def\VM{\Lfd{V}{M}}
\def\Xa{\Lfd{X}{1}}
\def\Xb{\Lfd{X}{2}}
\def\Xc{\Lfd{X}{3}}
\def\Xd{\Lfd{X}{4}}
\def\Xe{\Lfd{X}{5}}
\def\Yi{\Lfd{Y}{i}}
\def\Bb{\Bf{b}}
\def\Bg{\Bf{g}}
\def\Bu{\Bf{u}}
\def\Bust{\Bfu{u}{*}}
\def\Bui{\Bfd{u}{i}}
\def\Bv{\Bf{v}}
\def\Bvist{\Bfdu{v}{i}{*}}
\def\Bt{\Bf{t}}
\def\Bz{\Bf{z}}
\def\Bo{\Bf{0}}
\def\Blk{\Bfd{1}{\!k}}
\def\BA{\Bf{A}}
\def\BAlam{\Bfd{A}{\!\lam}}
\def\BAlamml{\Bfdu{A}{\!\lam}{-1}}
\def\BAml{\Bfu{A}{-1}}
\def\BB{\Bf{B}}
\def\BC{\Bf{C}}
\def\BD{\Bf{D}}
\def\BDdag{\Bfu{D}{\dagger}}

\def\BDRl{\Bfdu{D}{\xi}{\lam}}
\def\BDxi{\Bfd{D}{\xi}}
\def\BDxisq{\Bfdu{D}{\xi}{2}}
\def\BDximl{\Bfdu{D}{\xi}{-1}}
\def\BDa{\Bfd{D}{a}}
\def\BG{\Bf{G}}
\def\BHlam{\Bfd{H}{\!\lam}}
\def\BI{\Bf{I}}
\def\BIk{\Bfd{I}{\!k}}
\def\BIn{\Bfd{I}{\!n}}
\def\BIp{\Bfd{I}{\!p}}
\def\BP{\Bf{P}}

\def\BPpP{\BP'\!\BP}
\def\BQ{\Bf{Q}}
\def\BR{\Bf{R}}
\def\BU{\Bf{U}}
\def\BUst{\Bfu{U}{*}}
\def\BUa{\Bfd{U}{1}}
\def\BV{\Bf{V}}
\def\BVst{\Bfu{V}{*}}
\def\BW{\Bf{W}}
\def\BX{\Bf{X}}
\def\BXdag{\Bfu{X}{\dagger}}
\def\BXM{\Bfd{X}{\!\lam}}
\def\BXpX{\BX'\!\BX}
\def\BXpXml{(\BXpX)^{-1}}
\def\BXpXM{\BXM'\!\BXM}

\def\BY{\Bf{Y}}
\def\BYhat{\widehat{\BY}}
\def\BYhatc{\BYhat_{\!c}}
\def\BYhatRl{\BYhat_{\!R_{\lam}}}
\def\BYhato{\BYhat_{\!0}}

\def\BZ{\Bf{Z}}
\def\BZpZ{\BZ'\!\BZ}
\def\BZpZml{(\BZpZ)^{-1}}
\def\BZM{\Bfd{Z}{\!\lam}}
\def\BZpZM{\BZM'\!\BZM}
\def\alfbar{\bar{\alf}}
\def\alfa{\Lfd{\alf}{1}}

\def\alfk{\Lfd{\alf}{k}}

\def\bet{\Lf{\beta}}

\def\betasq{\Lfdu{\bet}{1}{2}}
\def\betksq{\Lfdu{\bet}{k}{2}}
\def\bo{\Lfd{\beta}{0}}
\def\ba{\Lfd{\beta}{1}}
\def\bb{\Lfd{\beta}{2}}
\def\bc{\Lfd{\beta}{3}}
\def\bd{\Lfd{\beta}{4}}
\def\be{\Lfd{\beta}{5}}
\def\bhat{\widehat{\beta}}
\def\bhatLa{\bhat_{L_{\!1}}}

\def\bhatLj{\bhat_{L_{\!j}}}
\def\bhatLk{\bhat_{L_{\!k}}}

\def\bhata{\bhat_{1}}
\def\bhatb{\bhat_{2}}
\def\bhatc{\bhat_{3}}
\def\bhatd{\bhat_{4}}
\def\bhate{\bhat_{5}}
\def\bhatk{\bhat_{k}}
\def\eps{\Lf{\epsilon}}
\def\epsi{\Lfd{\eps}{i}}
\def\lamst{\Lfu{\lam}{*}}
\def\lamhat{\widehat{\lam}}
\def\lamoc{\Lfd{\lam}{c}}
\def\lamocst{\Lfdu{\lam}{c}{*}}
\def\lama{\Lfd{\lam}{1}}
\def\lamb{\Lfd{\lam}{2}}
\def\lamk{\Lfd{\lam}{k}}

\def\muo{\Lfd{\mu}{0}}
\def\muSc{\Lfd{\mu}{S_{c}}}
\def\muBc{\Lfd{\mu}{\Bc}}

\def\omi{\Lfd{\om}{i}}

\def\phibu{\phi(\cdot)}
\def\psibu{\psi(\cdot)}
\def\phia{\Lfd{\phi}{1}}
\def\phio{\Lfd{\phi}{0}}
\def\phimin{\Lfd{\phi}{min}}
\def\phist{\Lfu{\phi}{*}}
\def\sigsq{\sigma^{2}}
\def\sighatsq{\widehat{\sigma}^{2}}
\def\xia{\Lfd{\xi}{1}}
\def\xie{\Lfd{\xi}{5}}
\def\xiasq{\Lfdu{\xi}{1}{2}}
\def\xii{\Lfd{\xi}{i}}
\def\xiisq{\Lfdu{\xi}{i}{2}}
\def\xik{\Lfd{\xi}{k}}
\def\xiksq{\Lfdu{\xi}{k}{2}}
\def\Blam{\Bf{\Lam}}
\def\bbet{\Bf{\beta}}
\def\bbeto{\Bfd{\beta}{\!0}}
\def\bbeta{\Bf{\beta}}
\def\bbhat{\widehat{\bbet}}
\def\bbhatRl{\bbhat_{\!R_{\lam}}}
\def\bbhatSl{\bbhat_{\!S_{\lam}}}
\def\bbhatc{\bbhat_{\!c}}
\def\bbhato{\bbhat_{\!0}}
\def\bbhatst{{\bbhat^{*}}}
\def\bbhatstl{\bbhat^{*}_{\!\lam}}
\def\bbhatLS{\bbhat_{\!L}}
\def\bbhatR{\bbhat_{\!R}}
\def\beps{\Bf{\eps}}
\def\bmu{\Bf{\mu}}
\def\bth{\Bf{\theta}}
\def\bthhat{\widehat{\bth}}
\def\bthhatLS{\bthhat_{\!L}}
\def\bthhatR{\bthhat_{\!R}}
\def\bthhatRl{\bthhat_{\!R_{\lam}}}
\def\bthhatSl{\bthhat_{\!S_{\lam}}}

\def\NbbhatRl{\parallel\nmsp\bbhatRl\nmsp\parallel}

\def\Blam{\Bf{\Lambda}}

\def\Bsig{\Bf{\Sigma}}

\begin{document}
\title[]{{\textbf{ANOMALIES IN THE FOUNDATIONS\\
OF RIDGE REGRESSION}}}
\author{D. R. Jensen and D. E. Ramirez}
\address{Department of Mathematics\\
         University of Virginia\\
         Charlottesville, VA 22904-4137}
\email{der@virginia.edu}
\thanks{}

\date{}
\begin{abstract}
Anomalies persist in the foundations of ridge regression as set
forth in Hoerl and Kennard (1970) and subsequently. Conventional
ridge estimators and their properties do not follow on constraining
lengths of solution vectors using LaGrange's method, as claimed.
Estimators so constrained have singular distributions; the proposed
solutions are not necessarily minimizing; and heretofore
undiscovered bounds are exhibited for the ridge parameter. None of
the considerable literature on estimation, prediction,
cross--validation, choice of ridge parameter, and related issues,
collectively known as \emph{ridge regression,} is consistent with
constrained optimization, nor with corresponding inequality
constraints. The problem is traced to a misapplication of LaGrange's
principle, failure to recognize the singularity of distributions,
and misplaced links between constraints and the ridge parameter.
Other principles, based on condition numbers, are seen to validate
both conventional ridge and \emph{surrogate ridge} regression to be
defined. Numerical studies illustrate that ridge analysis often
exhibits some of the same pathologies it is intended to redress.
\end{abstract}
\maketitle{}

\section{Introduction} \label{S:intro}

Given the full--rank model $\BY = \BX\bbet + \beps$ with zero--mean,
homoscedastic, and uncorrelated errors, the ordinary least squares
\emph{(OLS)} estimators $\bbhatLS$ solve the $k$ equations
$\BXpX\bbet = \BX'\BY$ on minimizing $Q(\bbet) =
(\BY-\BX\bbet)'(\BY-\BX\bbet).$ Ill--conditioned models long have
posed special challenges, in that $\bbhatLS$ often exhibits
excessive length, inflated variances, instability, and other
intrinsic difficulties. Noting these, Hoerl (1962, 1964) considered
\emph{ad hoc} solutions $\bbhatR$ = $\{\bbhatRl = (\BXpX +
\lam\BIk)^{-1}\BX'\BY; \lam\geq 0\}$ and noted their successful
applications in chemical engineering. Analyses built around these
have been labeled \emph{ridge regression} in statistics, although
Levenberg (1944) and Riley (1955) earlier posed such solutions in
numerical analysis. Noting that \emph{OLS} ``does not have built
into it a method for portraying sensitivity of the solutions to the
estimation criterion,'' Hoerl and Kennard (1970) sought mathematical
foundations beyond Gauss's principle with its inherent limitations.
Specifically, they asserted that $\bbhatR$ are solutions minimizing
$Q(\bbet)$ subject to the constraint $\{\bbet'\bbet = \csq\}.$
Others identify ridge regression instead with the constraints
$\{\bbet'\bbet \leq \csq\}$ of Balakrishnan (1963); however, Hoerl
and Kennard (1970), p. 64, specifically relegate this to approaches
other than ridge regression.

Ridge estimators abound, based on estimative, predictive,
cross--validative, and numerous other criteria, typically giving
disparate choices for $\lam.$ Even the early simulations of
Dempster, Schatzoff, and Wermuth (1977) identified 57 ridge and
related shrinkage estimators. An expository survey and numerical
examples are provided in Myers (1990). In short, a considerable
literature, spanning the past thirty--six years, rests on the
foundations of Hoerl and Kennard (1970), ostensibly the mathematics
of constrained optimization, to remedy defects of \emph{OLS} in
ill--conditioned systems.

In fact, little of the collective literature known as ridge
regression is consistent with the constrained optimization of Hoerl
and Kennard (1970), nor with corresponding inequality constraints.
Here the problem is traced to (i) a misapplication of LaGrange's
principle, (ii) failure to identify singular distributions, and
(iii) invalid links between the constraints and the ridge
parameters. These errors are evident also in Marquardt (1970),
Marquardt and Snee (1975), Golub, Heath and Wahba (1979), van Nostrand (1980), and
elsewhere throughout the literature. In consequence, much that is
known about ridge regression rests on a false premise. By analogy,
Hoerl and Kennard (1970) considered generalized ridge regression as
solving the modified equations $(\BXpX + \Blam)\bbet = \BX'\BY,$
with nonnegative ridge parameters $\Blam =
\textup{Diag}(\seq{\lam}{k}).$ As noted later, these solutions again
are inconsistent with LaGrange minimization. In summary, not to
denigrate its usefulness in practice, the collective body of ridge
regression rests on little more than heuristics. To the contrary,
aspects of ridge regression have proven useful enough, often enough,
to deserve sound rationale for their implementation. In this spirit
we seek to supplant the missing foundations with alternatives based
on conditioning of the linear system $\BXpX\bbet = \BX'\BY$. An
outline follows.

Supporting developments comprise Section 2, to include notation and
the basics of invariance and condition numbers. Section 3 reexamines
LaGrange optimization in linear inference. Section 4 develops
supporting rationale for ridge regression as currently practiced,
and an alternative approach using \emph{surrogate ridge} models. A
case study in Section 5 revisits an ill--conditioned data set
considered elsewhere. Section 6 concludes with a brief summary.

\section{Preliminaries}\label{S:prelim}

\subsection{Notation}\label{SS:nota}
The symbols $\Rk$ and $\Rkp$ designate Euclidean $k$--space and its
positive orthant; $\Fnk$ and $\FnkC$ comprise the real and complex
$(n\times k)$ matrices of rank $k\leq n;$ and $\Sk$ and $\Skp$
designate the real symmetric $(k\times k)$ matrices and their
positive definite varieties. The transpose, inverse, trace, and
determinant of $\BA\in\Fkk$ are $\BA',$ $\BAml,$ $\textup{tr}(\BA),$
and $\mid\negmedspace\BA\!\mid,$ and $\BVst$ is the conjugate
transpose of $\BV\in\FnkC.$ Groups of note include $\Uk$ as the
unitary $(k\times k)$ matrices, and $\Ok$ as the real orthogonal
group. Special arrays are the $(k\times k)$ identity $\BIk,$ the
unit vector $\Blk$ = $[1,1,\,\dots,\,1]'\in\Rk,$ and the diagonal
matrix $\BDa$ = $\textup{Diag}(\seq{a}{k}).$ The mapping $\sig(\BX)$
= $[\seq{\xi}{k}]'$ takes $\BX\in\FnkC$ into its ordered singular
values $\{\xia\!\geq\!\ldots\!\geq\!\xik\!>0\}.$ The
\textit{singular decomposition} is $\BX = \BU\BD\BVst,$ such that
$\BD$ = $\textup{Diag}(\BDxi,\Bo)$ of order $(n\times k),$ $\BDxi$ =
$\textup{Diag}(\seq{\xi}{k}),$ $\BU\in\Un,$ and $\BV\in\Uk,$ where
the columns of $\BU$ = $[\seq{\Bu}{n}]$ and $\BV$ = $[\seq{\Bv}{k}]$
comprise the \textit{left--} and \textit{right--singular vectors} of
$\BX.$ Equivalently write $\BX = \BUa\BDxi\BVst$ =
$\Sumilk\xii\Bui\Bvist$ with $\BUa$ = $[\seq{\Bu}{k}],$ and its
Moore--Penrose inverse as $\BXdag$ = $\BV\BDdag\BUst,$ with $\BDdag$
= $\textup{Diag}(\BDxi^{-1},\Bo)$ of order $(k\times n).$
Specifically, the real model $\BY = \BX\bbet + \beps$ in canonical
form becomes $\BY = \BP\BDxi\bth + \beps,$ where $\BX =
\BP\BDxi\BQ'$ and $\bth = \BQ'\bbet$ is an orthogonal
reparametrization. For $\BY\in\Rn$ random, designate its mean
vector, its dispersion and correlation matrices as $\textup{E}(\BY)
= \bmu,$ $\textup{V}(\BY) = \Bsig,$ and $\BC(\BY)$ = $\BR,$ and its
law of distribution as $\La(\BY).$

\subsection{Invariance and Conditioning.}\label{SS:invarcond}
A function $\psibu$ on $\FnkC$ is called \textit{unitarily
invariant} if, for each $\BG\in \FnkC$ and any unitary matrices
$\BU\in\Un$ and $\BV\in\Uk,$ it follows that $\psi(\BG)$ =
$\psi(\BU\BG\BVst).$ Then $\psi(\BG)$ depends on $\BG$ only through
its ordered singular values $\sig(\BG) = [\seq{\gam}{k}]'.$ Let
$\Phi$ comprise the \textit{symmetric gauge functions} on $\Rk$ such
that for each $\phibu\in\Phi,$ (i) $\phi(\seq{u}{k})$ is symmetric
under the $2^{k}k!$ permutations and reflections about the origin;
(ii) $\phi(\Bu) > 0$ when $\Bu \neq \Bo;$ (iii) $\phibu$ is
homogeneous, \emph{i.e.,} $\phi(c\Bu)$ = $\mid\!\!c\!\!\mid\phi(\Bu
)$ for $c \neq 0;$ and (iv) $\phi(\Bu + \Bv ) \leq \phi(\Bu) +
\phi(\Bv).$ Let $\Psi$ comprise the unitarily invariant matrix norms
on $\FnkC;$ von Neumann (1937) demonstrated that these are generated
as $\{\parallel\!\cdot\!\parallel_{\phi};\, \phi\in\Phi\}$ with
$\parallel\!\BG\!\parallel_{\phi}$ = $\phi(\seq{\gam}{k}).$
Corresponding norms on $\Fnk$ are invariant under
$\BX\to\BP\BX\BQ',$ with $(\BP,\,\BQ)\in\On\!\times\!\Ok;$ see also
Schatten (1970) and Marshall and Olkin (1979). In particular, the
Frobenius norm on $\Fnk$ is $\parallel\!\BX\!\parallel_{\!F}$ =
$[\textup{tr}(\BXpX)]^{1/2}$ = $(\Sumilk\xiisq)^{1/2}$ in terms of
the singular decomposition $\BX = \BP\BDxi\BQ',$ with
$\parallel\!\cdot\!\parallel$ as the Euclidean norm on $\Rk.$

Two types of conditioning are germane to the present study:

\textit{Type A Conditioning:} Stability of the solution $\Bz$ of the
linear system $\BA\Bz = \Bb,$ when the coefficients $\BA\in\Fkk$ are
subjected to small perturbations, is gauged by the \textit{condition
number} $\cg(\BA)$ = $g(\BA)g(\BAml),$ where $\gbu$ ordinarily is a
norm. The system is well conditioned at $\BA$ = $\BIk$ with
$\cg(\BIk)$ = 1.0, larger values reflecting greater
ill-conditioning. Specifically, with $g(\BA)$ =
$\parallel\!\BA\!\parallel_{\phi},$ then $\{\cphibu;\,\phi\in\Phi\}$
comprise the \emph{unitarily invariant Type A condition numbers,} so
that $\{\cphi(\BA)$ =
$\parallel\!\!\BA\!\parallel_{\phi}\parallel\!\!\BAml\!\!
\parallel_{\phi};\,\phi\in\Phi\},$
as treated in Marshall and Olkin (1979), Horn and Johnson (1985),
and elsewhere. In particular, take $\ca(\BA) = \alfa/\alfk,$ where
$\{\alfa\!\geq\ldots\geq\!\alfk\}$ are the ordered eigenvalues of
$\BA.$

\textit{Type B Conditioning:} The concept of elasticities is invoked
in Belsley, Kuh and Welsch (1980) to link sensitivities of
solutions, and of variances of $\bbhatLS$ = $\BZpZml\BZ'\BY,$ with
disturbances in the data matrix $\BZ\in\Fnk,$ as gauged by its
condition number $\ca(\BZ)$ = $\xia/\xik,$ with $\sig(\BZ) =
[\seq{\xi}{k}]'.$ Here $\BZ$ is the result of scaling the columns of
$\BX$ to have (approximately) equal lengths. More generally, the
unitarily invariant condition numbers on $\Fnk$ are $\cphi(\BX)$ =
$\phi(\BX)\phi(\BXdag),$ with $\BXdag$ as the Moore--Penrose
inverse. The system is well conditioned at $\BX$ = $\BP\BIk\BQ',$
where $\cphi(\BP\BIk\BQ')$ = $\cphi(\BIk)$ = 1.0, larger values
reflecting greater ill-conditioning. In summary, Belsley \textit{et
al.} (1980) proceed to scale the columns of $\BX\to\BZ$ to have
approximately equal lengths, and to focus on
$\parallel\!\!\BZ\!\!\parallel_{\phia}$ = $\xia,$ so that $\ca(\BZ)$
= $\xia/\xik.$

\newpage
\section{The Principal Issues} \label{S:pissues}

\subsection{LaGrange's Method.}\label{SS:LaGMeth}
Given differentiable functions $f(\seq{u}{k})$ and $g(\seq{u}{k})$
such that the gradient $\nabla g(\seq{u}{k})\neq\Bo$ on $\Go =
\{\Bu\in\Rk: g(\Bu) = 0\},$ the problem is to minimize
$f(\seq{u}{k})$ subject to the constraint $g(\seq{u}{k}) = 0$. Write
$L(\seq{u}{k},\lam)$ = $f(\seq{u}{k})$ + $\lam[g(\seq{u}{k})-0].$ It
is necessary that gradient vectors in $\Rk$ be parallel,
\textit{i.e,}
\begin{equation}
\nabla f(\seq{u}{k}) = \lam\nabla g(\seq{u}{k}),
\end{equation}
\noindent whereas
\begin{equation}
\partial L(\seq{u}{k},\lam)/\partial \lam = [g(\seq{u}{k}) - 0]
\end{equation}
\noindent recovers the constraint. LaGrange's principle requires
solving the $k+1$ equations, (3.1) and (3.2), in the $k+1$ unknowns
$\{\seq{u}{k},\lam\}.$ To minimize $f(\seq{u}{k})$ subject to
$g(\seq{u}{k})\geq 0,$ define the Lagrangian $L(\Bu,\lam)$ = $f(\Bu)
- \lam g(\Bu).$ Stuezle\footnote{Stuezle, W., ``Chapter 5. Notes on
Ridge Regression,'' online notes for BioStat 538, Winter 2005,
University of Washington, at
www.stat.washington.edu/wxs/Stat538--w05} has given conditions for
$\Bust$ to be a solution, namely, (i) $g(\Bust)\geq 0;$ (ii)
$\nabla_{\Bu} L(\Bust,\lamst) = \Bo;$ (iii) $\lamst g(\Bust) = \Bo;$
and (iv) $\lamst \geq 0.$

For constrained least squares the objective function is now
\begin{equation}
L(\seq{\beta}{k},\lam) = Q(\bbeta) + \lam(\bbet'\bbet - \csq)\notag
\end{equation}
\noindent with $Q(\bbeta)$ = $(\BY-\BX\bbet)'(\BY-\BX\bbet)$ as
before. Corresponding to (3.1) and (3.2) are
\begin{equation}
(\BXpX + \lam\BIk)\bbet = \BX'\BY
\end{equation}
\begin{equation}
\bbet'\bbet = \csq
\end{equation}
\noindent to be solved for the $k+1$ unknowns
$(\seq{\beta}{k},\lam).$ Designate these as $\{\bbhatc,\lamhat\}$
such that $\bbhatc'\bbhatc = \csq,$ as apparently intended by Hoerl
and Kennard (1970). If instead $Q(\bbeta)$ is to be minimized
subject to $\{\bbet'\bbet\leq\csq\},$ then the constrained solution
$\bbhato$ satisfies $\bbhato = \bbhatLS$ whenever $\bbhatLS'\bbhatLS
< \csq,$ and otherwise $(\BXpX + \lam\BIk)\bbeto = \BX'\BY$ for some
$\lam > 0$ such that $\bbhato'\bbhato = \csq,$ as shown in
Balakrishnan (1963). See also the conditions (i)--(iv) of Stuezle
(2005) as cited.

\subsection{Ridge Regression: A Survey.}\label{SS:RRSur}
We recall essentials of conventional ridge regression as set forth
principally in Hoerl and Kennard (1970), Marquardt (1970), and
Marquardt and Snee (1975). For continuity we retain their notation,
with their $\{k,p,\bbhat, \bbhatst\}$ corresponding to our $\{\lam,
k, \bbhatLS, \bbhatRl\}$ and, on occasion, $\bbhatRl$ = $\bbhatstl$
= $\bbhatst\!(\lam).$ Accordingly, write the residual sum of squares
as $\phi$  = $(\BY-\BX\Bb)'(\BY-\BX\Bb)$ = $\phi_{min}$ +
$\phi(\Bb),$ where $\phimin$ = $(\BY-\BX\bbhat)'(\BY-\BX\bbhat)$ and
$\phi(\Bb)$ = $(\Bb-\bbhat)'\BXpX(\Bb-\bbhat).$ Various assertions
have been set forth, as enumerated here for later reference.

\textbf{A1.} Hoerl and Kennard (1970), p. 57: ``$\bbhatst = [\BIp +
k\BXpXml]^{-1}\bbhat$ (2.3)*.''

\textbf{A2.} Hoerl and Kennard (1970), pp. 58--59: ``The ridge trace
can be shown to be following a path through the sums of squares
surface so that for a fixed $\phi$ a single value for $\Bb$ is
chosen and that is the one with minimal length.'' Precisely:
``Minimize $\Bb'\Bb$ subject to $(\Bb-\bbhat)'\BXpX (\Bb-\bbhat) =
\phio$ (3.2)*.'' ``This reduces to $\Bb$ = $\bbhatst$ =
$(\BXpX+k\BI)^{-1}\BX'\BY$ where $k$ is chosen to satisfy the
restraint (3.2)*.''

\textbf{A3.} Hoerl and Kennard (1970), p. 59: ``Of course, in
practice it is easier to choose a $k\geq 0$ and then to compute
$\phio.$ In terms of $\bbhatst$ the residual sum of squares becomes
$\phist(k$) = $(\BY-\BX\bbhatst)'(\BY-\BX\bbhatst)$ = $\phimin +
\ksq\bbhatst{'}\BXpXml\bbhatst$ (3.6)*.''

\textbf{A4.} Hoerl and Kennard (1970), p. 59: ``A completely
equivalent statement of the problem is this: If the squared length
of the regression vector $\Bb$ is fixed at $\Rsq,$ then $\bbhatst$
is the value of $\Bb$ that gives a minimum sum of squares. That is,
$\bbhatst$ is the value of $\Bb$ that minimizes the function $\Fa$ =
$(\BY-\BX\Bb)'(\BY-\BX\Bb) + (1/k)(\Bb'\Bb - \Rsq)$ (3.7)*.''

\textbf{A5.} Marquardt and Snee (1975), p. 5: ``If $\bbhatst$ is the
solution of $(\BXpX + k\BI)\bbhatst = \Bg,$ then $\bbhatst$
minimizes the sum of squares of residuals on the sphere centered at
the origin whose radius is the length of $\bbhatst.$'' Here
$\Bg=\BX'\BY.$

\subsection{Properties of Solutions.}\label{SS:PropsSol}
We next examine distributions of the constrained solutions
$\bbhatc,$ subject to $\bbhatc'\bbhatc = \csq,$ to continue the
unfinished work of Hoerl and Kennard (1970), and of $\bbhato$ under
inequality constraints. To these ends identify the sphere $\Sc =
\{\Bu\in\Rk: \Bu'\Bu = \csq\}$ and the open ball $\Bc = \{\Bu\in\Rk:
\Bu'\Bu < \csq\},$ both of radius $c,$ and the complement $\Bcc =
\{\Bu\in\Rk: \Bu'\Bu \geq \csq\}.$ Accordingly, let $\mu(\cdot)$ be
the probability measure on $\Rk$ induced through $\BY\to\bbhatLS;$
let $\muSc(\cdot)$ be the measure on $\Sc\subset\Rk$ induced through
solutions $\bbhatc$ of (3.3) and (3.4); and let $\muBc(\cdot)$ be
the nonsingular measure on $\Bc\subset\Rk$ induced through
$\La(\bbhatLS\negmedspace\mid\negmedspace \bbhatLS'\bbhatLS\!<\!
\csq).$ Stochastic properties of $\bbhatc$ and $\bbhato$ are given
next.

\begin{theorem}
Let $\bbhatc\in\Rk$ be the constrained solution satisfying (3.3) and
(3.4), and let $\bbhato\in\Rk$ minimize $Q(\seq{\beta}{k})$ subject
to $\{\bbet'\bbet\leq\csq\},$ with $\muo(\cdot)$ as its probability
measure on $\Rk.$

(i) The  joint distribution $\La(\bbhatc)$ = $\Fc(\Bb)$
corresponding to $\muSc(\cdot)$ is singular on $\Rk$ of rank $k-1.$

(ii) The measure $\muo(\cdot)$ for $\bbhato$ admits the mixture
representation
\begin{equation}
\muo(A) = \alf\cdot\muBc(A) + \alfbar\cdot\muSc(A)
\end{equation}
\noindent with mixing probabilities $\alf = 1 - \alfbar\in (0,1),$
such that

(iii) $\muBc(A) = [\mu(\Bc)]^{-1}\int_{A} I_{\Bc}(\Bt)d\mu(\Bt),$
where $I_{\Bc}(\Bt)$ is the indicator function; and

(iv) $\alf = \mu(\Bc).$
\end{theorem}
\noindent\textbf{Proof:} Conclusion (i) is immediate, since
$\bbhatc\in\Rk$ constructively lies on the sphere $\bbhatc'\bbhatc =
\csq.$ We proceed by conditioning on the exclusive outcomes
$\bbhatLS\in\Bc$ and $\bbhatLS\in\Bcc.$ Clearly $\bbhato$ takes the
value $\bbhatLS$ with probability $\alf = P(\bbhatLS'\bbhatLS\!<\!
\csq)$ = $\mu(\Bc),$ where the conditional measure corresponding to
$\La(\bbhatLS\negmedspace\mid\negmedspace \bbhatLS'\bbhatLS\!<\!
\csq)$ is $\muBc(A)$ = $[\mu(\Bc)]^{-1}\int_{A}
I_{\Bc}(\Bt)d\mu(\Bt),$ as asserted, to give conclusion (iii).
Similarly, $\bbhato$ takes the value $\bbhatc$ with probability
$\alfbar = 1-\alf$ as in (iv), its conditional measure as in (i), to
complete our proof. $\square$

Observe that the singular distribution $\La(\bbhatc)$ of conclusion
(i) may be added to the list of distributions arising in the
analysis of directional data, to include the von Mises--Fisher
distributions, for example. For further reference see Batschelet
(1981), Fisher (1993), Fisher, Lewis and Embleton (1993), Evans,
Hastings and Peacock (2000), and Mardia and Jupp (2000). Conclusion
(ii) for $\bbhato$ complements the work of Balakrishnan (1963) in
the context of linear estimation. Moreover, under Gaussian errors,
$\alf = \mu(\Bc)$ derives from a weighted sum of $k$ independent
noncentral chi--squared random variables, each having a single
degree of freedom; see Kotz, Johnson and Boyd (1967).

\subsection{A Critique.}\label{SS:Critiq}
We next reexamine the critical assertions of Section 3.2.

Assertion \textbf{A1:} False.  As noted, the solution $\bbhatc$
necessarily lies on the sphere $\bbhatc'\bbhatc = \csq$ and thus has
a joint singular distribution in $\Rk$ of rank $k-1.$ To the
contrary, Assertion \textbf{A1} implies that $\bbhatst(\lam)$ has a
nonsingular distribution for each $\lam\geq 0$, yet $\bbhatst$
clearly refers to the constrained solution throughout Section 3 of
Hoerl and Kennard (1970). The assertion is false, applying to
solutions of (3.3) only, as there is no one--to--one linear
transformation taking $\bbhatLS$ onto the sphere $\bbhat'\bbhat =
\csq.$ In consequence, expression (3.6)* of Hoerl and Kennard (1970)
is in error, as are its implications, since the term
$\ksq\bbhatst{'}\BXpXml\bbhatst$ derives from the inapplicable
Assertion \textbf{A1.}

Assertion \textbf{A2,} and its dual \textbf{A4,} appear to be
essentially intact. The exception is that ``$1/k$'' in expression
(3.7)* of Hoerl and Kennard (1970) instead should be ``$k.$''

Assertion \textbf{A5:} False. This assertion arises as the dual to
\textbf{A3,} excluding (3.6)* of Hoerl and Kennard (1970). The basic
idea is to solve (3.3) as $\bbhatst\!(\lam)$ for fixed $\lam > 0,$
and then to discover the implied constraint $\{\bbet'\bbet =
\cstsq\}$ at (3.4) on evaluating $\bbhatst{'}\bbhatst = \cstsq.$
However, the solution $\bbhatst\!(\lam)$ need not minimize the
residual sum of squares $\SS(\lam)$ =
$[\BY-\BX\bbhatst\!(\lam)]'[\BY-\BX\bbhatst\!(\lam)],$ as claimed.
This fallacy stems from the tacit but unfounded assumption that
$\lam$ and $\csq$ correspond one--to--one. To the contrary, it is
demonstrated in Section 5 that multiple solutions may have the same
length but different $\lam$s, for example,
$\parallel\!\bbhatst\!(\lama)\!\parallel$ =
$\parallel\!\bbhatst\!(\lamb)\!\parallel$ with $\lama < \lamb.$ But
then the solution $\bbhatst\!(\lamb)$ cannot be minimizing, as
$SS(\lamb) > SS(\lama)$ from the monotonicity of $SS(\lam).$ In this
regard Figure 3 of Marquardt and Snee (1975) is particularly
misleading. Assertions \textbf{A2,} ``for a fixed $\phi$ a single
value for $\Bb$ is chosen and that is the one with minimal length,''
and \textbf{A5,} that ``$\bbhatst$ minimizes the sum of squares of
residuals on the sphere centered at the origin whose radius is the
length of $\bbhatst,$'' often are misrepresented as equivalent
assertions regarding solutions $\bbhatRl$ of (3.3) alone. See van
Nostrand (1980), for example.

To continue, for fixed $c$ define the equivalence class
\begin{equation}
\Lam(c) = \{\lam :\, \parallel\!\bbhatst\!(\lam)\!\parallel = c\},
\end{equation}
\noindent and let $\lamoc$ = $\min\{\Lam(c)\}.$ Then Assertion
\textbf{A5} may be corrected as follows.

Assertion \textbf{A5*.} If $\bbhatst\!(\lam)$ is a solution of
$(\BXpX + \lam\BI)\bbhatst = \BX'\BY$ having length
$\parallel\!\bbhatst\!(\lam)\!\parallel = \cst,$ then
$\bbhatst(\lamocst)$ minimizes the sum of squares of residuals on
the sphere centered at the origin whose radius is the length $\cst$
of $\bbhatst,$ where $\lamocst$ = $\min\{\Lam(\cst)\}.$

Assertion \textbf{A5*} has profound consequences in practice. Of the
many schemes devised for choosing the ridge parameter $\lam,$ the
user then must examine the corresponding equivalence class for each
such $\lam.$ If it is a singleton set, then the solution thus
attained is minimizing. Otherwise the algorithm \textbf{A5*} must be
implemented to attain the minimizing solution. Further details are
provided in Section 5.3.

It is clear that $\bbhatc$ is the LaGrange solution minimizing
$Q(\bbeta)$ subject to $\{\bbet'\bbet = \csq\}.$ To the contrary,
Hoerl and Kennard (1970), Marquardt (1970), Marquardt and Snee
(1975), Golub \emph{et al.} (1979), and others concerned with
constrained optimization, instead take $\bbhatRl$ as the ridge
estimator, solving (3.3) alone for some $\lam > 0.$ Together with
Assertion \textbf{A5,} this is tantamount to asserting that the $k$
linear equations (3.3) somehow embody the constraint (3.4) as well,
which they clearly cannot. Yet $\bbhatR,$ not $\bbhatc,$ comprise
\emph{the} ridge estimators on which essentially all of ridge
regression now rests. Assertion \textbf{A1} clearly holds for
solutions $\bbhatRl$ satisfying (3.3) only.

Confusion persists in the meaning of ridge regression. Bunke (1975),
Hocking (1976), and Tibshirani (1996), for example, assert that
ridge regression embodies the inequality constraint
$\{\bbet'\bbet\leq\csq\},$ despite the disclaimer of Hoerl and
Kennard (1970). Yet nowhere do these authors acknowledge the
constrained solution $\bbhato$ of Balakrishnan (1963), nor its
properties as in Theorem 1, opting instead for the ridge solutions
$\{\bbhatRl; \lam\geq 0\}$ of Hoerl (1962, 1964). On the other hand,
the inequality--constrained solution $\bbhato$ does have the
nonsingular mixture distribution of Theorem 1. However, we are aware
of no work in ridge regression that explicitly accounts for the
structure of either $\bbhatc$ or of $\bbhato$ as in Theorem 1.

In short, ridge regression in its present form rests essentially on
$\bbhatR$ through an accident of history. Indeed, expressions for
variances and biases; solutions for $\lam$ purporting to minimize
expected mean squares; prediction, validation, and
cross--validation; and other aspects of ridge regression; all are
predicated on Assertion \textbf{A1.} If instead either $\bbhatc$ or
$\bbhato$ were taken as starting points, as required under the aegis
of constrained optimization, then the ensuing ``ridge regressions''
would differ dramatically from the conventional one based on
$\{\bbhatRl; \lam\geq 0\},$ together with the critical but false
Assertion \textbf{A1,} with $\csq$ now corresponding to $\lam.$
These differences necessarily would include issues such as (i) the
stability of the solutions $\bbhatc$ or $\bbhato$ instead of
$\bbhatRl,$ in comparison with $\bbhatLS;$ (ii) the inflation of
variances, taking into account actual variances to be derived from
Theorem 1 as reference; (iii) prediction using $\BYhatc =
\BX\bbhatc$ or $\BYhato = \BX\bbhato,$ instead of $\BYhatRl =
\BX\bbhatRl;$ (iv) the use, meaning, and properties of
cross--validative and predictive criteria based on $\BYhatc$ or
$\BYhato,$ instead of $\BYhatRl;$  (v) ridge traces as modified to
take into account $\bbhato$ and singularity of the joint
distribution of $\bbhatc;$ and (vi) the trade--off between bias and
variance of the constrained estimators $\bbhatc$ and $\bbhato,$ as
determined using actual moments to be derived from Theorem 1. Other
differences may be noted. All such properties would have to be
established anew, complicated considerably by the nonstandard
distributions encountered in Theorem 1.

By analogy, Hoerl and Kennard (1970) further considered generalized
ridge regression invoking the $k$ equations $(\BXpX + \Blam)\bbet =
\BX'\BY,$ with $\Blam = \textup{Diag}(\seq{\lam}{k})$ as nonnegative
ridge parameters. Note that this, too, cannot have resulted from
LaGrange minimization: Given that
$\{\betasq=\casq,\ldots,\betksq=\cksq\},$ the only function of the
data now would be to determine signs of the roots $\{\bhata = \pm
\ca,\ldots,\bhatk = \pm\ck\}.$ On the other hand, if inequality
constraints $\{\betasq\leq\casq,\ldots,\betksq\leq\cksq\}$ are
invoked instead, then correct solutions are provided by Myoken and
Uchida (1977) akin to those of Balakrishnan (1963) where $\{\lama =
\cdots = \lamk = \lam\}.$

\section{Foundations Via Conditioning}\label{S:cndfound}

We seek substitutes for the failed principle of constrained
optimization as a basis for conventional ridge regression. In what
follows we consider $\{\bbhatRl;\lam\geq 0\}$ as solutions to (3.3)
alone as in Hoerl (1962, 1964), without reference to constrained
optimization and discredited assertions thereto as noted. Type A
conditioning of the linear system $\BXpX\bbet = \BX'\BY$ prompts the
modification $(\BXpX + \lam\BIk)\bbet = \BX'\BY,$ from the
perspective of both numerical analysis (Levenberg (1944) and Riley
(1955)) and of statistics (Hoerl (1962, 1964)). A survey is provided
subsequently. Moreover, the Type B conditioning of $\BY =\BX\bbet +
\beps$ is also germane, since the conditioning of $\BXpX$ depends on
that of $\BX,$ and for further reasons to be cited. A new approach
to ill conditioned systems, using \textit{surrogate ridge models,}
rests essentially on Type B conditioning. Details follow.

\subsection{Background}\label{SS:Bkgd}
Ill--conditioned models typically arise from nonorthogonality of
columns of $\BX.$ Let $\BW = \BXpX$ and $\BV = \BXpXml.$ Since
$\textup{V}(\bbhatLS)$ = $\sigsq\BV,$ the \emph{variance inflation
factors} (\emph{VIF}s) of $\bbhatLS$ = $[\bhatLa,\ldots,\bhatLk]'$
are defined as $\{V\!I\!F(\bhatLj) = \vjj/\wjj^{-1}; 1\leq j\leq
k\},$ \emph{i.e.,} the ratio of the actual variance to the ``ideal''
variance attained when columns of $\BX$ are orthogonal, so that
$\BW$ = $\textup{Diag}(\waa,\ldots,\wkk).$ Often $\BY = \BZ\bbet +
\beps$ is taken with $\BZpZ$ in ``correlation form'' having unit
diagonal elements; then $\{V\!I\!F(\bhatLj) = \vjj; 1\!\leq\!
j\leq\! k\}$ are diagonal elements of $\BV = \BZpZml$ from the
scale--invariance of \emph{VIF}s. With
$\{\Va\geq\Vb\geq\ldots\geq\Vk\}$ as the \emph{ordered} diagonal
elements of $\BV,$ Marquardt and Snee (1975) identify $\Va$ to be
``the best single measure of the conditioning of the data,'' thus a
critical diagnostic tool. See also Marquardt (1970), Beaton, Rubin
and Barone (1976), and Davies and Hutton (1975).  A basic connection
between \emph{VIF}s and condition numbers is due to Berk (1977):
\begin{lemma}
Given $\BZpZ$ in correlation form, with
$\{\Va\!\geq\!\Vb\!\geq\!\ldots\!\geq\!\Vk\}$ as the ordered
diagonal elements of $\BV = \BZpZml.$ Then the condition number
$\ca(\BZpZ)$ satisfies
\begin{equation}
\Va \leq \ca(\BZpZ) \leq k(\Va + \cdots + \Vk).
\end{equation}
\end{lemma}
\noindent Since $\{\cphi(\BA) = \cphi(\BAml); \phi\in\Phi\}$ from
Section 2.2, the Type A condition number for $\BZpZ\bbet = \BZ'\BY$
is identical to $\cphi[\textup{V}(\bbhatLS)],$ so that Lemma 1 is
really about dispersion parameters in the equivalent form
\begin{equation}
\Va \leq \ca[\textup{V}(\bbhatLS)] \leq k(\Va + \cdots + \Vk).
\end{equation}

\subsection{Ridge Regression}\label{SS:RR}
That $\BXpX \to (\BXpX + \lam\BIk)$ improves conditioning has been
cited by Marshall and Olkin (1979) as a justification for ridge
regression. In brief, their Theorem C.3, p. 273, asserts that for
any $(\BA,\BB)\in\Skp$ such that $\cphi(\BB)\leq\cphi(\BA),$ with
$\{\cphi(\cdot);\phi\in\Phi\}$ as in Section 2.2, then $\cphi(\BA +
\BB)\leq\cphi(\BA).$ Riley (1955) showed that $\BB = \lam\BIk$
satisfies the hypothesis of the theorem for any $\BA\in\Skp,$ where
$\lam$ depends on numerical considerations. This holds for any Type
A conditioning of $\BA\Bz = \Bb$ $\to$ $(\BA + \lam\BIk)\Bz = \Bb$
as in Section 2.2, and thus in particular for $\BXpX\bbet = \BX'\BY$
$\to$ $(\BXpX + \lam\BIk)\bbet = \BX'\BY,$ as noted by Marshall and
Olkin (1979), p. 273, to give Type A conditioning as a basis for
ridge regression. Moreover, using condition numbers $\ca(\cdot),$
the improvement is seen directly on comparing $\ca(\BXpX)$ =
$\xiasq/\xiksq$ with $\ca(\BXpX + \lam\BIk)$ = $(\xiasq +
\lam)/(\xiksq + \lam),$ where $\sig(\BX) = [\seq{\xi}{k}]'.$
Essential properties of $\bbhatLS$ and $\bbhatRl$ are summarized in
Table 1, along with the \emph{surrogate estimator,} $\bbhatSl,$ to
be defined subsequently.

\begin{table}[ht]
\caption{Properties of $\{\bbhatLS,\bbhatRl, \bbhatSl\}$ under
Gauss--Markov assumptions, where $\BXM$ =
$\BP\textup{Diag}(\sqrt{\xiasq+\lam},\ldots,\sqrt{\xiksq+\lam})\BQ'$
and $\BAlam = (\BXpX + \lam\BIk).$}
\begin{center}
\begin{tabular}{| c | c | c | c |}\hline
Estimator & Definition & $\textup{E}(\bbhat)$ &
$\textup{V}(\bbhat)$\\
\hline
 $\bbhatLS$ & $\BXpXml\BX'\BY$ & $\bbet$ & $\sigsq\BXpXml$\\
 $\bbhatRl$ & $\BAlamml\BX'\BY$ & $\BAlamml\BX'\BX\bbet$ &
$\sigsq\BAlamml\BX'\BX\BAlamml$\\
 $\bbhatSl$ & $\BAlamml\BXM'\BY$ & $\BAlamml\BXM'\BX\bbet$&
 $\sigsq\BAlamml$\\
\hline
\end{tabular}
\end{center}
\end{table}

\subsection{Surrogate Models}\label{SS:SR}
Nonetheless, the correspondence $\BA\Bz = \Bb$ $\longleftrightarrow$
$\BXpX\bbet = \BX'\BY$ is incomplete in the context of linear
inference, since both $\BA = \BXpX$ and $\Bb = \BX'\BY$ are subject
to disturbances in $\BX.$ This has not been taken into account. In
particular, ridge solutions satisfying $(\BXpX + \lam\BIk)\bbet =
\BX'\BY,$ despite improved conditioning on the left, still are
subject to the ill--conditioning of $\BX$ on the right. To correct
this oversight, we invoke Type B conditioning from Section 2.2 on
observing that $\BXpX\to(\BXpX + \lam\BIk)$ is tantamount to
modifying $\BX$ itself as a means to enhanced conditioning. In
particular, begin with the singular decomposition $\BX =
\BP\BDxi\BQ';$ let $\BXM$ =
$\BP\textup{Diag}(\sqrt{\xiasq+\lam},\ldots,\sqrt{\xiksq+\lam})\BQ';$
observe that $(\BXpX + \lam\BIk)$ = $\BXM'\BXM;$ and note that ridge
regression entails $\BXM'\BXM\bbet = \BX'\BY.$ Instead, we take $\BY
= \BXM\bbet + \beps$ as an approximation, or \textit{surrogate,} for
the ill--conditioned model $\BY = \BX\bbet + \beps$ itself, as in
the following.

\begin{definition}
Given an ill--conditioned model $\BY = \BX\bbet + \beps,$  its
\emph{ridge surrogate} is a modified model $\BY = \BXM\bbet +
\beps.$ The \emph{surrogate estimator} $\bbhatSl,$ solving
$\BXM'\BXM\bbet = \BXM'\BY,$ is \emph{OLS} for the surrogate model.
\end{definition}

To continue, the order of approximation of $\BXM$ for $\BX$ may be
gauged by the Frobenius distance
\begin{equation}
\parallel\!\BX\!-\!\BXM\!\parallel_{\!F} = \left[\Sumilk\left(\xii -
\sqrt{\xiisq+\lam}\right)^{2}\right]^{1/2},
\end{equation}
\noindent from the unitary invariance of
$\parallel\!\cdot\!\parallel_{\!F}.$ Moreover, the conditioning of
$\BXM'\BXM\bbet = \BXM'\BY$ now may be gauged through Type B
conditioning as in Section 2.2. For later reference, basic
properties of $\{\bbhatSl; \lam\geq 0\}$ are summarized in Table 1.
It remains to compare properties of $\{\bbhatLS, \bbhatRl,
\bbhatSl\}.$ Direct comparisons are somewhat obscure; however, these
become more transparent on invoking canonical forms to be considered
next.

\subsection{Canonical Forms}\label{SS:CF}
The singular decomposition $\BX = \BP\BDxi\BQ',$ with $\BPpP =
\BIk,$ together with the orthogonal reparametrization $\bth =
\BQ'\bbet,$ gives $\BY = \BX\bbet + \beps$ $\to$ $\BY =
\BP\BDxi\BQ'\bbet + \beps$ $\to$ $\BU = \BP'\BY = \BDxi\bth +
\BP'\beps,$ such that $\textup{E}(\BP'\beps) = \Bo$ and
$\textup{V}(\BP'\beps)$ = $\sigsq\BP'\BIn\BP$ = $\sigsq\BIk$ under
Gauss--Markov assumptions regarding the errors of $\BY = \BX\bbet +
\beps.$ Accordingly, $\textup{E}(\BU) = \BDxi\bth$ and
$\textup{V}(\BU)$ = $\sigsq\BIk.$ In canonical form it follows that
$\bthhatLS = (\BDxisq)^{-1}\BDxi\BU$ = $\BDximl\BU,$
$\textup{E}(\bthhatLS) = \bth,$ and $\textup{V}(\bthhatLS)$ =
$\sigsq\BDxi^{-2}$ under \emph{OLS,} as given in Table 2. Similar
expressions for the canonical ridge estimators $\{\bthhatRl;
\lam\geq 0\},$ and the canonical surrogate ridge estimators
$\{\bthhatSl; \lam\geq 0\},$ are reported in Table 2. Since $\bbhat
= \BQ\bthhat,$ $\textup{E}(\bbhat) = \BQ\textup{E}(\bthhat),$ and
$\textup{V}(\bbhat) = \BQ\textup{V}(\bthhat)\BQ'$ for all three
estimators, Table 1 follows directly from Table 2, and conversely.
Moreover, issues regarding the conditioning of $\{\bbhatLS,
\bbhatRl, \bbhatSl\},$ as linear data transformations, and
conditioning of the corresponding dispersion matrices
$\{\textup{V}(\bbhatLS), \textup{V}(\bbhatRl),
\textup{V}(\bbhatSl)\},$ are considered subsequently. These can be
established directly in terms of those of $\{\bthhatLS, \bthhatRl,
\bthhatSl\},$ since $\BQ$ is orthogonal and condition numbers here
are unitarily invariant.

\begin{table}[ht]
\begin{center}
\caption{Properties of $\{\bthhatLS, \bthhatRl, \bthhatSl\}$ under
standard Gauss--Markov assumptions, where $\BU = \BP'\BY$ and
$\BD(\omi) = \textup{Diag}(\seq{\om}{k}).$}
\renewcommand{\arraystretch}{1.25}
\begin{tabular}{| c | c | c | c |}\hline
Estimator & Definition & $\textup{E}(\bthhat)$ &
$\textup{V}(\bthhat)$\\
\hline
 $\bthhatLS$ & $\BDxi^{-1}\BU$ & $\bth$ & $\sigsq\BDxi^{-2}$\\
 $\bthhatRl$ & $\BD(\xii/(\xiisq+\lam))\BU$ & $\BD(\xiisq/(\xiisq+\lam))\bth$ &
$\sigsq\BD(\xiisq/(\xiisq+\lam)^{2})$\\
 $\bthhatSl$ & $\BD(1/\sqrt{\xiisq+\lam})\BU$ & $\BD(\xii/(\sqrt{\xiisq+\lam}))\bth$&
 $\sigsq\BD(1/(\xiisq+\lam))$\\
\hline
\end{tabular}
\end{center}
\end{table}

Specifically, in canonical form we have $\BDxi\bthhatLS = \BU,$ so
that the Type B condition number $\ca(\BX)$ = $\ca(\BDxi)$ =
$\xia/\xik$ properly gauges the sensitivity of the solution
$\bbhatLS$ to disturbances in $\BX.$ Similarly, with $\BDRl$ =
$\textup{Diag}((\xiasq+\lam)/\xiasq,\ldots,(\xiksq+\lam)/\xiksq),$
observe from $\BDRl\bthhatRl = \BU$ that its condition number gauges
sensitivity of the solution $\bthhatRl,$ and thus of $\bbhatRl =
\BQ\bthhatRl$ to perturbations in $\BX,$ from the orthogonality of
$\BQ.$ This underscores the central role of Type B conditioning from
Section 2.2, as set forth in Belsley \emph{et al.} (1970).

\subsection{Central Issues}\label{SS:CI}
Several issues, to be examined empirically in Section 5, appear to
be open questions not addressed in the voluminous literature on
ridge regression. Intrinsic difficulties with \emph{OLS} include (i)
nonorthogonality of the columns of $\BX,$ as reflected in
$\cphi(\BX)$ and $\cphi(\BXpX);$ (ii) instability of solutions
linked to the conditioning of the data transformation
$\bbhatLS(\BY)$ = $\BXpXml\BX'\BY,$ considered as a function of
$\BY;$ and (iii) pathologies in dispersion parameters as reflected
in \emph{VIF}s and the ill--conditioning of $\textup{V}(\bbhatLS).$
Moreover, at some level the conditioning of $\textup{E}(\bbhatRl)$ =
$T(\bbet)$ becomes an issue in transforming the parameter space, as
in assessing the trade-off between variance and bias. As ridge
regression seeks remedies, it is pertinent to ask how well the ridge
solutions progress towards those ends. Regarding item (i), the
apparent ``correlations'' in $\BW = \BXpX,$ namely
$\{\wij/\sqrt{\wii\wjj}\},$ are taken into
$\{\wij/\sqrt{(\wii+\lam)(\wjj+\lam)}\}$ as elements of $(\BXpX +
\lam\BIk).$ These in turn decrease in magnitude with increasing
$\lam.$ Nonetheless, ridge solutions themselves are subject to
nonorthogonality, together with attendant difficulties regarding
stability, \emph{VIF}s, and conditioning of their dispersion
matrices. Improving stability of the solutions thus hinges on the
conditioning of $\bbhatRl(\BY)$ when considered as a data
transformation. Moreover, the capacity to ameliorate dispersion
problems of \emph{OLS} hinges on improving \emph{VIF}s and condition
numbers for $\textup{V}(\bbhatRl).$ On the other hand, it is widely
known that $\bbhatRl$ shrinks stochastically towards the origin, as
do its mean and dispersion matrix, with increasing $\lam$. These
issues in turn prompt several questions to be considered
subsequently.

\begin{description}
\item[Q1] Does it follow that stability of $\bbhatRl(\BY)$ necessarily
improves with increasing $\lam?$
\item[Q2] Given that $\textup{V}(\bbhatRl)$ = $\sigsq(\BXpX +
\lam\BIk)^{-1}\BXpX(\BXpX + \lam\BIk)^{-1},$ does it follow that
condition numbers $\ca[\textup{V}(\bbhatRl)]$ decrease with
increasing $\lam?$
\item[Q3] With regard to variance inflation, does it follow that \emph{VIF}s
for elements of $\bbhatRl$ decrease with increasing $\lam?$
\item[Q4] Viewing $\textup{E}(\bbhatRl)$ =
$T(\bbet)$ as a transformation on the space of parameters, does it
follow that its conditioning improves with increasing $\lam?$
\end{description}

For completeness, observe that the foregoing issues pertain not only
to the ridge estimators $\{\bbhatRl; \lam\geq 0\}$ themselves, but
also to other biased solutions to include $\{\bbhatSl; \lam\geq
0\}$.

We next undertake a comparative study of properties of ridge and
surrogate ridge solutions, to be continued in the case studies of
Section 5.

\subsection{Some Comparisons}\label{SS:CF}
Regarding the conventional $\{\bbhatRl; \lam\geq 0\}$ and surrogate
ridge $\{\bbhatSl; \lam\geq 0\}$ estimators, both shrink
stochastically towards the origin with increasing $\lam,$ as do
their means and variances, and similarly for $\{\bthhatRl; \lam\geq
0\}$ and $\{\bthhatSl; \lam\geq 0\}.$ Specifically, for a given
$\lam,$ it is seen  from Table 2 that $\bthhatSl$ achieves lesser
shrinkage, both in expectation and variance, than $\bthhatRl.$

Condition numbers for various arrays are given in Table 3 for the
canonical estimators $\{\bthhatLS, \bthhatRl, \bthhatSl\}.$ These
arrays include (i) coefficients defining $\bthhat(\BU)$ with
reference to stability of the solutions; (ii) coefficients defining
the parameter transformations $\textup{E}(\bthhat) = T(\bth);$ and
(iii) the dispersion matrix $\textup{V}(\bthhat).$ Entries in Table
3 follow directly from Table 2 and the definition of $\ca(\cdot),$
on recalling that elements of $\BDxi$ =
$\textup{Diag}(\seq{\xi}{k})$ are ordered as
$\{\xia\geq\ldots\geq\xik > 0\}.$ Observe, moreover, that the rows
of Table 3 may be identified equivalently as $\{\bbhatLS, \bbhatRl,
\bbhatSl\},$ and the columns as $\{\ca[\bbhat(\BY)], \ca[T(\bbet)],
\ca[\textup{V}(\bbhat)]\},$ respectively. This follows since $\bbet
= \BQ\bth,$ $\bbhat = \BQ\bthhat,$ and $\textup{V}(\bbhat)$ =
$\BQ\textup{V}(\bthhat)\BQ',$ $\BQ$ is orthogonal, and the condition
numbers are unitarily invariant.

\begin{table}[ht]
\begin{center}
\caption{Condition numbers for data transformations $\bthhat(\BU),$
for parameter transformations $\textup{E}(\bthhat)$ = $T(\bth),$ and
for $\textup{V}(\bthhat),$ for each of $\{\bthhatLS, \bthhatRl,
\bthhatSl\}.$}
\renewcommand{\arraystretch}{1.5}
\begin{tabular}{| c | c | c | c |}\hline
Estimator & $\ca[\bthhat(\BU)]$ & $\ca[T(\bth)]$ &
$\ca[\textup{V}(\bthhat)]$\\
\hline
 $\bthhatLS$ & $\frac{\xia}{\xik}$ & 1.00 & $\frac{\xiasq}{\xiksq}$\\
 $\bthhatRl$ &
 $\frac{\max\{\xii/(\xiisq+\lam)\}}{\min\{\xii/(\xiisq+\lam)\}}$ &
 $\frac{\xiasq(\xiksq+\lam)}{\xiksq(\xiasq+\lam)}$ &
 $\frac{\max\{\xiisq/(\xiisq+\lam)^{2}\}}{\min\{\xiisq/(\xiisq+\lam)^{2}\}}$\\
 $\bthhatSl$ & $\frac{\sqrt{\xiasq+\lam}}{\sqrt{\xiksq+\lam}}$ &
$\frac{\xia\sqrt{\xiksq+\lam}}{\xik\sqrt{\xiasq+\lam}}$
& $\frac{\xiasq+\lam}{\xiksq+\lam}$\\
\hline
\end{tabular}
\end{center}
\end{table}

Note further that $\ca[\bbhatLS(\BY)]$ = $\ca(\BX)$ and
$\ca[\textup{V}(\bbhatLS)]$ = $\ca(\BXpX),$ whereas
$\ca[\bbhatSl(\BY)]$ = $\ca(\BXM)$ and $\ca[\textup{V}(\bbhatSl)]$ =
$\ca(\BXpXM),$ as both are \emph{OLS} in their respective models.
Moreover, both condition numbers, $\ca[\bbhatSl(\BY)]$ =
$(\sqrt{\xiasq+\lam}/\sqrt{\xiksq+\lam}),$ and its square
$\ca[\textup{V}(\bbhatSl)],$ decrease monotonically with increasing
$\lam,$ thus assuring improved conditioning for the surrogate
estimators. Condition numbers associated with $\bthhatRl,$ and thus
with $\bbhatRl,$ are more convoluted and will be examined further in
Section 5.

\section{Case Studies} \label{S:casestud}
\subsection{The Data}\label{SS:data}
We reexamine the Hospital Manpower Data as reported in Myers (1990).
Records at $n=17$ U. S. Naval Hospitals include: $Y:$ Monthly
man--hours; $\Xa:$ Average daily patient load; $\Xb:$ Monthly X--ray
exposures; $\Xc:$ Monthly occupied bed days; $\Xd:$ Eligible
population in the area $\div$ 1000; and $\Xe:$ Average length of
patients' stay in days. The basic model is
\begin{equation}
\Yi = \bo + \ba\Xa + \bb\Xb + \bc\Xc + \bd\Xd + \be\Xe + \epsi;
1\leq i \leq n.
\end{equation}
\noindent Following Hoerl and Kennard (1970), Marquardt (1970),
Marquardt and Snee (1975), Myers (1990), and others, we center and
scale the model, so that $\BY = \BZ\bbet + \beps$ with $\BZpZ$ in
correlation form, the central focus being the rates of change
$\bbet$ = $[\ba, \bb, \bc, \bd, \be]'.$ The data are given in Table
3.8, pp. 132--133, of Myers (1990), and computations were done
mostly using PROC IML of the SAS Programming System. The data are
exceedingly ill--conditioned: Elements of $\BDxi$ are $\BDxi$ =
$\textup{Diag}(2.048687,\, 0.816997,\, 0.307625,\, 0.201771,\,
0.007347);$ $\ca(\BZpZ)$ = 77,754.86; the maximal \emph{VIF} in
\emph{OLS} estimation is $\Va$ = $V\!I\!F(\bhata)$ = 9,595.685; and
other \emph{VIF}s appear subsequently in Table 8 at $\lam = 0.$

\subsection{Choices for $\lam$}\label{SS:CStud}
Widely diverse criteria have evolved in the choice for $\lam,$ with
profound consequences regarding ridge estimators, ridge predictors,
and their properties. Five criteria in common usage are reported in
Table 4,
\begin{table}[ht]
\begin{center}
\caption{Choices for $\lam$ in the Hospital Manpower Data
corresponding to conventional criteria $\DFlam,$ $\GCVlam,$ $\Clam,$
$\PRElam,$ and $\HKBlam.$}
\renewcommand{\arraystretch}{1.3}
\begin{tabular}{| c | c | c |}\hline
 Name & \hspace{1.1in}Definition\hspace{1.1in} & Value for $\lam$\\ \hline
 $\DFlam$ & $\textup{tr}(\BHlam)$ = $\Sumilk\frac{\xiisq}{(\xiisq+\lam)}$ & 0.0004\\
 $\GCVlam$ & $\frac{SS_{Res,\lam}}{[n-(1+\textup{tr}(\BHlam))]^{2}}$ & 0.004787\\
 $\Clam$ & $[\frac{SS_{Res,\lam}}{\sighatsq}-n+2+2\textup{tr}(\BHlam)]$ & 0.0050\\
 $\PRElam$ & $\Sumiln\eilsq$ & 0.2300\\
 $\HKBlam$  &  $\frac{k\sighatsq}{\bbhatLS'\bbhatLS}$ & 0.616964\\ \hline
\end{tabular}
\end{center}
\end{table}
together with definitions and their values as determined for the
Hospital Manpower Data. These include $\DFlam = \textup{tr}(\BHlam)$
with $\BHlam$ = $[\BZ(\BZpZ + \lam\BIk)^{-1}\BZ'];$ the
cross--validation $\PRElam$ statistic of Allen (1974); a
rotation--invariant version called Generalized Cross Validation
($\GCVlam$) by Golub \emph{et al.} (1979); $\Clam$ as a device for
variance--bias trade--off as in Mallows (1973); and $\HKBlam$ as
recommended by Hoerl, Kennard and Baldwin (1975) based on simulation
studies. As listed in Table 4, $SS_{Res,\lam}$ is the residual sum
of squares using ridge regression; $\sighatsq$ is the \emph{OLS}
residual mean square; and $\{\eilsq\}$ are the $P\!R\!E\!S\!S$
residuals for ridge regression. Further details are given in Myers
(1990), pp. 392--411, including numerical values for $\DFlam,$
$\Clam,$ and $\PRElam$ as reported in Table 4. Further choices
include $\lam\in\{0.01, 0.03, 0.05, 0.07, 0.09\}$ and others to be
noted subsequently.

\subsection{Minimizing Solutions}\label{SS:MinSol}
Often a definitive value for the constraint $\{\bbet'\bbet = \csq\}$
is not apparent in a particular study. This motivates the dual
Assertions \textbf{A3} and \textbf{A5} of Section 3.2: (i) choose
$\lam;$ (ii) solve (3.3) for $\bbhatRl;$ (iii) evaluate the implied
constraint at (3.4) as $\bbhatRl'\bbhatRl = \cstsq;$ and (iv) assert
as in \textbf{A5} that the solution so attained ``minimizes the sum
of squares of residuals on the sphere centered at the origin whose
radius is the length'' of $\bbhatRl.$ We have claimed that Assertion
\textbf{A5} is false. Evidence is provided in Table 5,
\begin{table}[ht]
\begin{center}
\caption{Lengths of $\bbhatRl,$ and square roots of residual sums of
squares $\Rlam =
[(\BY-\BZ\bbhatRl)'(\BY-\BZ\bbhatRl)]^{\frac{1}{2}},$ for designated
values of $\lam.$}
\begin{tabular}{| c | c | c | c | c | c | c | c | c |}\hline
 $\lam$ & 0.00 & 0.04 & 0.08 & 0.12 & 0.16 & 0.20 & 0.24 & 0.28\\ \hline
 $\NbbhatRl$ & 394.67 & 137.82 & 33.14 & 31.50 & 70.02 & 99.19 & 122.10 & 140.73\\
 $\Rlam$ & 2129.53 & 2474.87 & 2735.75 & 2914.38 & 3057.54 & 3184.84 & 3305.00 & 3422.13\\ \hline
 $\lam$ & 0.32 & 0.36  & 0.40 & 0.48 & 0.56 & 0.60 & 0.64 & 0.68\\ \hline
 $\NbbhatRl$ & 156.25 & 169.40 & 180.69 & 199.00 & 213.09 & 218.93 & 224.11 & 228.70\\
 $\Rlam$ & 3538.22 & 3654.22 & 3770.58 & 4004.70 & 4240.27 & 4358.28 & 4476.26 & 4594.06\\ \hline
 $\lam$ & 0.72 & 0.76 & 0.80 & 0.84 & 0.88 & 0.92 & 0.96 & 1.00\\  \hline
 $\NbbhatRl$ & 232.79 & 236.42 & 239.65 & 242.53 & 245.08 & 247.34 & 249.33 & 251.09\\
 $\Rlam$ & 4711.56 & 4828.65 & 4945.23 & 5061.20 & 5176.49 & 5291.03 & 5404.77 & 5517.64\\ \hline
\end{tabular}
\end{center}
\end{table}
where lengths $\NbbhatRl,$ and square roots $\Rlam$ = $[(\BY -
\BZ\bbhatRl)'(\BY - \BZ\bbhatRl)]^{\frac{1}{2}},$ are reported as
$\lam$ ranges systematically over $[0, 1].$ Recall that this range
is stipulated by Hoerl and Kennard (1970) and others when $\BZpZ$ is
in ``correlation form.'' Here $\bbhatRl$ =
$[\bhata,\bhatb,\bhatc,\bhatd,\bhate]'$ consists of rates of change;
similar trends are exhibited when $\bbet$ is expanded to include the
intercept. It is seen that $\NbbhatRl$ initially decreases to a
minimum, then increases beyond $\lam$ = 1.0, but eventually
decreases to zero since $\bbhatRl$ is a shrinkage estimator.

Greater detail is seen on recalling from Section 4.4 that $\bbhatRl
= \BQ\bthhatRl;$ that $\BQ$ is orthogonal; and thus, letting
$g_{\bbhatR}(\lam)$ = $\NbbhatRl^{2},$ that $g_{\bbhatR}(\lam)$ =
$g_{\bthhatR}(\lam).$ The canonical form of Section 4.4 assures that
$g_{\bthhatR}(\lam)$ = $\Sumilk\Uisq\xiisq/(\xiisq + \lam)^{2}.$
This is differentiable; its derivative is
\begin{equation}
\partial g_{\bthhatR}(\lam)/\partial \lam = -2\Sumilk\Uisq
\xiisq(\xiisq + \lam)^{-3};
\end{equation}
\noindent and its path traces evolution of the derivative as $\lam$
varies. In particular, at $\lam = 0$ we have $[\partial
g_{\bthhatR}(\lam)/\partial\lam]_{\!\lam = 0}$ =
$-2\Sumilk\Uisq/\xii^{4}.$  This is precipitous for the Hospital
Manpower Data in view of the fact that $\xik = \xie$ = 0.007347.

A detailed local view is provided in Table 6, to include not only
$\NbbhatRl$ and $\Rlam,$ but also the ridge estimates $\bbhatRl$ =
$[\bhata,\bhatb,\bhatc,\bhatd,\bhate]'$ in rows corresponding to
various choices for $\lam.$
\begin{table}[ht]
\begin{center}
\caption{Ridge estimators $\bbhatRl,$ lengths of $\bbhatRl,$ and
square roots $\Rlam =
[(\BY-\BZ\bbhatRl)'(\BY-\BZ\bbhatRl)]^{\frac{1}{2}}$ of residual
sums of squares, for designated values of $\lam.$}
\begin{tabular}{| c | c | c | c | c | c | c | c |}\hline
$\lam$ & $\bhata$ & $\bhatb$ & $\bhatc$ & $\bhatd$ & $\bhate$ &
$\NbbhatRl$ & $\Rlam$\\  \hline
 0.08 & 10.6354 & 0.065428 & 0.359139 & 6.3206 & -30.7471 & 33.1448 & 2735.75\\
 0.08095 & 10.6118 & 0.065432 & 0.358279 & 6.3674 & -28.9649 & 31.5000 & 2740.68\\
 0.08797 & 10.4475 & 0.065444 & 0.352298 & 6.6903 & -16.4728 & 20.6250 & 2775.83\\
 0.0981 & 10.2378 & 0.065414 & 0.344681 & 7.0942 & -0.3156 & 12.4645 & 2823.03\\
 0.09829 & 10.2342 & 0.065413 & 0.344548 & 7.1012 & -0.0308 & 12.4615 & 2823.89\\
 0.0983 & 10.2340 & 0.065413 & 0.344541 & 7.1015 & -0.0159 & 12.4615 & 2823.93\\
 0.11 & 10.0248 & 0.065325 & 0.336955 & 7.4935 & 16.3900 & 20.6251 & 2874.22\\
 0.12 & 9.8679 & 0.065217 & 0.331280 & 7.7785 & 28.8834 & 31.5000 & 2914.38\\
\hline
\end{tabular}
\end{center}
\end{table}
Values of $\bbhatRl$ for $\lam\in\{0.08,\, 0.11,\, 0.12\}$ are as in
Table 8.9 of Myers (1990), who reports ridge estimates for $\lam\in
[0,\, 0.24]$ by increments of 0.01. It is seen that $\NbbhatRl$
takes its minimum value, 12.46150, at $\lam_{\min}$ = 0.09829. To
continue, designate $\bbhatRl$ as $\bbhatR(\lam).$ It is seen that
$\bbhatR(0.12)$ and $\bbhatR(0.08095)$ have the same length, namely,
$\parallel\nmsp\bbhatR(0.12)\nmsp\parallel$ = 31.500 =
$\parallel\nmsp\bbhatR(0.08095)\nmsp\parallel,$ so that
$\Lam(31.500) = \{0.08095,\, 0.12\}$ in the notation of (3.6).
Suppose a user chooses $\bbhatR(0.12)$ as the ridge estimate for the
Hospital Manpower Data. Then $\bbhatR(0.12)$ is not the minimizing
solution of length 31.500; this is seen from $R(0.12) = 2914.38 >
2740.68 = R(0.08095).$ Similarly, it is clear that $\Lam(20.625) =
\{0.08797,\, 0.11\}$ as in (3.6), and that $\bbhatR(0.11)$ of Table
6 is not minimizing, to be supplanted instead by $\bbhatR(0.08797)$
from Table 6. A continuum of further examples can be constructed by
reflecting $\lam$ asymmetrically about $\lam_{\min}$ = 0.09829, the
smaller $\lam$ of each pair corresponding to the minimizing
solution. These clearly constitute counterexamples to Assertion
\textbf{A5.}

Not only are definitive values for the constraint $\{\bbet'\bbet =
\csq\}$ not evident beforehand, but profound and heretofore
undiscovered limits pertain to admissible values for $\lam$ in order
that solutions of given length $\cst$ be minimizing. To fix ideas,
suppose in equation (3.4) that $\{g_{\bbhatR}(0.00) >
\cstsq\geq\csq\geq g_{\bbhatR}(0.09829) = 12.46150^{2} =
155.2877\}.$ Then the only feasible values for $\lam$ are those in
the interval $[\min g_{\bbhatR}^{-1}(\csq),\,0.09829].$ For example,
if $33.14481^{2} = 1098.5784\geq\csq\geq 155.2877,$ then from Table
6 the feasible values are $\lam\in[0.08,\,0.09829].$ For
$\{g_{\bbhatR}(0.00)\geq\csq\geq g_{\bbhatR}(0.09829)\},$ the
feasible values are $\lam\in[0.00,\,0.09829].$ These are the only
feasible values  for $\lam\in[0,\,1].$ On the other hand, choosing
$\{0 < \csq < g_{\bbhatR}(0.09829) = 155.2877\}$ requires $\lam$ in
the interval $(\min g_{\bbhatR}^{-1}(\csq),\,\infty),$ where $\min
g_{\bbhatR}^{-1}(155.2877) > 158.$ For example, if $\csq < 100,$
then the feasible values are $\lam\in (198,\,\infty).$ As these are
far outside the recommended interval $[0,\,1],$ constraints
$\cstsq\in (0,\,155.2877)$ must be declared to be inadmissible.
Values reported for $\{0 < \csq < g_{\bbhatR}(0.09829) = 155.2877\}$
are supported by the Maple software package. Values reported for
$\PRElam$ and $\HKBlam$ in Table 4 are thus inadmissible in view of
Assertion \textbf{A5*.}

In short, imbedded in the Hospital Manpower Data are the hidden
feasible constraints $\{\bbet'\bbet = \csq\}$ with $\csq \geq
155.2877.$ These could not have been discerned beforehand short of
the foregoing detailed analyses.

To summarize, origins of the anomaly exhibited here may be traced as
follows: (i) The ridge trace of $\bhate(\lam)$ exhibits a
down-up--down character, beginning with $\bhate(0.00)$ = -394.3280,
decreasing to zero between $\lam$ = 0.09 and $\lam$ = 0.10, and
increasing thereafter to $\bhate(1.00)$ = 250.8307 and beyond, and
eventually decreasing to zero through shrinkage. (ii)
$\mid\nmsp\bhate(\lam)\nmsp\mid$ dominates other estimates by orders
of magnitude ranging from one to four except near its minimum. (iii)
Other estimates exhibit relatively narrow ranges in comparison with
$\bhate(\lam)$ as $\lam$ varies over $[0,\,1].$ (iv) In consequence,
$\NbbhatRl$ is largely determined by $[\bhate(\lam)]^{2}$ as $\lam$
varies. Finally note that $\Lam(\cst)$ from (3.6) takes on two
values in the cases examined, from the down--up--down character of
$\NbbhatRl$ as $\lam$ evolves. It is clear in other circumstances
that $\Lam(\cst)$ may consist of three or more elements. For
example, a single dominant estimate may exhibit multiple sign
changes, whereas estimates for other coefficients may have one or
more sign changes as well. These and related matters are studied in
Zhang and McDonald (2005), and references cited therein, under
special structure of $\BZpZ$ in correlation form. Properties, to
include sign changes, crossings, and rates--of--change of individual
ridge estimates, as well as bounds on the number of sign changes,
are determined by those authors on identifying zeros and derivatives
of polynomials in $\lam$ of degree $k-1,$ under special structure as
cited.

These facts alone challenge the meaning of numerous simulation
studies purporting to compare alternative criteria for choosing
$\lam,$ when all such choices have ignored the minimizing
constraints on $\lam.$ Thus aggregates of minimizing/non--minimizing
values are compared with other such aggregates, to the effect of
total obfuscation.

We turn next to properties of ridge and surrogate ridge solutions,
to include condition numbers and other diagnostics. Computations for
the condition numbers proceed as in Table 3, based on equivalence
between conditioning for $\bbhat$ and the canonical estimators
$\bthhat,$ as noted in Section 4.6.

\subsection{Properties of $\bbhatRl$ and $\bbhatSl$}\label{SS:Props}
In summary, the ridge solutions to $(\BZpZ + \BIk)\bbet = \BZ'\BY$
account for ill--conditioning of $\BZpZ$ on the left of $\BZpZ\bbet
= \BZ'\BY,$ whereas the surrogate solutions to $\BZpZM\bbet$ =
$\BZM'\BY$ account for ill--conditioning on the right as well. It
thus is germane to compare $\{\bbhatSl;\lam\geq 0\}$ with
$\{\bbhatRl;\lam\geq 0\}$ using the data at hand. We next examine
 critical issues from Section 4.5, applicable both to ridge and
to surrogate ridge solutions. Table 7 lists condition numbers and
other quantities affiliated with $\{\bbhatRl;\lam\geq 0\}$ and
$\{\bbhatSl;\lam\geq 0\},$ under values for $\lam$ as listed.
\begin{table}[ht]
\begin{center}
\caption{Condition numbers for $\bbhatRl\!(\BY),$ $\bbhatSl\!(\BY),$
$\textup{V}(\bbhatRl),$ and $\textup{V}(\bbhatSl);$ the maximal
\emph{VIF}s $\VM(\bbhatRl)$ and $\VM(\bbhatSl);$ and the Frobenius
distance $\DZ(\BZM)$ = $\parallel\!\BZ-\BZM\!\parallel_{\!F},$ under
various choices for $\lam.$}
\begin{tabular}{| c | c | c | c | c | c | c | c |}\hline
$\lam$ &  $\ca(\bbhatRl)$  & $\ca(\bbhatSl)$ & $\VM(\bbhatRl)$ &
$\ca[\textup{V}(\bbhatRl)]$ & $\VM(\bbhatSl)$ &
$\ca[\textup{V}(\bbhatSl)]$ &
$\DZ(\BZM)$\\
\hline
 0.0004   &   33.1584 &   96.1565 & 141.5345 & 1099.4770 & 1146.399 & 9246.064  & 0.0140\\
 0.004787 &    9.0957 &   29.4630 &  10.9688 &   82.7319 & 112.6300 &  868.0653 & 0.0638\\
 0.005    &    9.0537 &   28.8348 &  10.8874 &   81.9695 & 108.0918 &  831.4473 & 0.0654\\
 0.010     &    8.1707 &   20.4561 &   9.2481 &   66.7610 &  56.6915 &  418.4530 & 0.0974\\
 0.030     &   11.6724 &   11.8596 &  21.2905 &  136.2440 &  21.2197 &  140.6508 & 0.1847\\
 0.050     &   15.1539 &    9.2114 &  34.0995 &  229.6392 &  13.6552 &   84.8507 & 0.2511\\
 0.070     &   17.8166 &    7.8046 &  42.5990 &  317.4320 &  10.2639 &   60.9119 & 0.3083\\
 0.090     &   20.4222 &    6.8997 &  51.7827 &  417.0673 &   8.3166 &   47.6061 & 0.3598\\
 0.230     &   29.6720 &    4.3868 & 100.5675 &  880.4276 &   3.9338 &   19.2438 & 0.6429\\
 0.616964 &   53.4183 &    2.7932 & 250.4309 & 2853.5130 &   2.0374 &    7.8022 & 1.1769\\
 1.000      &   66.6915 &    2.2797 & 451.5788 & 4447.7550 &   1.5976 &    5.1968 & 1.5745\\
\hline
\end{tabular}
\end{center}
\end{table}
Question 1 of Section 4.5 is negated for $\bbhatRl:$ Stability of
the solutions $\bbhatRl,$ as gauged by $\ca[\bbhatRl\!(\BY)],$
initially improves but then erodes. Further computations show that
$\ca[\bbhatRl\!(\BY)]$ takes its minimal value, 7.4463, at $\lam =
0.015,$ and increases thereafter. In contrast, despite higher
beginning values than $\ca[\bbhatRl\!(\BY)],$ the condition numbers
$\ca[\bbhatSl\!(\BY)]$ for surrogate estimators decrease
monotonically with increasing $\lam,$ the trends
$\ca[\bbhatRl(\BY)]$ = 11.7723 = $\ca[\bbhatSl(\BY)]$ crossing at
$\lam$ = 0.03045.

Questions 2 and 3 of Section 4.5 are refuted for $\bbhatRl:$
Computations interpolating those of Table 7 show that
$\ca[\textup{V}(\bbhatRl)]$ temporarily decreases over $\lam\in[0,\,
0.015],$ where its minimum is 55.4470, but it increases thereafter.
Similarly, the maximal \emph{VIF}s for $\bbhatRl$ initially decrease
and then increase. By comparison, both the condition numbers
$\ca[\textup{V}(\bbhatSl)],$ and the maximal \emph{VIF}s for
$\bbhatSl,$ decrease with increasing $\lam.$ Although initially
larger, $\VM(\bbhatSl)$ approximates $\VM(\bbhatRl)$ at $\lam =
0.030,$ and the ratio $\VM(\bbhatRl)/\VM(\bbhatSl)$ increases
markedly thereafter.

Recall that the surrogate $\BY = \BZM\bbet + \beps$ is intended as
an approximation to $\BY = \BZ\bbet + \beps.$ The order of
approximation, as gauged by the Frobenius distance (4.3), is
tabulated as the final column of Table 7. Relative changes, given by
$\parallel\!\BZ-\BZM\!\parallel_{\!F}\negmedspace/\negmedspace
\parallel\!\BZ\!\parallel_{\!F},$
are 0.1123 at $\lam = 0.05,$ ranging up to 0.5263 at $\lam =
0.616964,$ where the denominator is
$\parallel\!\BZ\!\parallel_{\!F}$ = 2.236068.

Further details are given in Tables 8 and 9, from which several
entries of Table 7 are drawn.
\begin{table}[ht]
\begin{center}
\caption{Variance inflation factors for $\bbhatRl,$ and condition
numbers for $\BC(\bbhatRl)$ and $T(\bbet)$ =
$\textup{E}(\bbhatRl),$ for designated values of $\lam.$}
\begin{tabular}{| c | c | c | c | c | c | c | c |}\hline
$\lam$ & \emph{VIF1} & \emph{VIF2} & \emph{VIF3} & \emph{VIF4} &
\emph{VIF5} & $\ca[\BC(\bbhatRl)]$ & $\ca[T(\bbet)]$ \\ \hline
    0.000    & 9595.68 & 7.9406 &  8931.449 & 23.2887 & 4.2794 & 54756.83 & 1.0000\\
    0.0004  & 141.5345 & 7.8481 & 133.0221 & 13.0512 & 3.3997 & 576.8409 & 8.4095\\
    0.004787  & 7.1604 & 7.1682   & 7.8840 & 10.9688 & 3.0128 & 90.13222 & 89.5726\\
    0.005    & 7.1047 & 7.1379 &    7.8349 & 10.8874 & 2.9972 & 89.50392 & 93.5175\\
    0.010  &    8.0001 & 6.4919 &    8.8456 & 9.2481 & 2.6830 & 75.66936 & 185.8150\\
    0.030  &   19.7743 & 4.7268 &   21.2905 & 5.6339 & 2.0003 & 109.4703 & 552.8219\\
    0.050  &   32.0013 & 3.6885 &   34.0995 & 4.0168 & 1.6988 & 177.2545 & 916.3722\\
    0.070  &   42.5990 & 3.0187 &   45.0695 & 3.1473 & 1.5363 & 227.9178 & 1276.515\\
    0.090  &   51.7827 & 2.5598 &   54.4589 & 2.6269 & 1.4377 & 267.3171 & 1633.297\\
    0.230  &  100.5675 & 1.3868 &  102.4723 & 1.6364 & 1.2446 & 441.9639 & 4040.511\\
    0.616964 & 250.4309 & 1.0879 & 243.2535 & 1.8791 & 1.3541 & 1047.931 & 9965.795\\
    1.000  &  451.5788 & 1.2184 &  430.5957 & 2.4738 & 1.5386 & 2174.418 & 14961.96\\
\hline
\end{tabular}
\end{center}
\end{table}
\begin{table}[ht]
\begin{center}
\caption{Variance inflation factors for $\bbhatSl,$ and condition
numbers for $\BC(\bbhatSl)$ and $\textup{V}(\bbhatSl),$ for
designated values of $\lam.$}
\begin{tabular}{| c | c | c | c | c | c | c |}\hline
$\lam$ & \emph{VIF1} & \emph{VIF2} & \emph{VIF3} & \emph{VIF4} &
\emph{VIF5} & $\ca[\BC(\bbhatSl)]$\\ \hline
 0.000    & 9595.68 & 7.9406 &  8931.449 & 23.2887 & 4.2794 & 54756.83\\
 0.0004 & 1146.399 &  7.8846 & 1068.211 & 14.2203 & 3.5089 & 5091.248\\
 0.004787 & 112.6300& 7.5308 & 106.0234 & 12.0987 & 3.2190 & 458.9380\\
 0.005 & 108.0918 &   7.5147 & 101.7946 & 12.0488 & 3.2099 & 440.5738\\
 0.010 & 56.6915 &      7.1607 & 53.8459 & 11.0379 & 3.0219 & 233.7461\\
 0.030 & 21.2197 &      6.0737 & 20.5412 & 8.4506  & 2.5374 & 93.1862\\
 0.050 & 13.6552 &      5.3181 & 13.3380 & 6.9511  & 2.2606 & 63.4167\\
 0.070 & 10.2639 &      4.7584 & 10.0777 & 5.9605  & 2.0792 & 48.8365\\
 0.090 & 8.3166 &        4.3258 & 8.1934 & 5.2538  & 1.9500 & 39.9124\\
 0.230 & 3.9338 &        2.8218 & 3.9092 & 3.1133  & 1.5493 & 17.9102\\
 0.616964&2.0374 &      1.7669 & 2.0330 & 1.8380  & 1.2710 & 7.5371\\
 1.000 & 1.5976 &         1.4635 & 1.5957 & 1.4981  & 1.1781 & 5.0614\\
\hline
\end{tabular}
\end{center}
\end{table}
Table 8 examines the evolution of \emph{VIF}s, and conditioning of
the correlation matrices, for $\bbhatRl$ as $\lam$ varies. Values
for $\ca[\BC(\bbhatRl)]$ are included, as Lemma 1 applies in each
case. It is found that $\ca[\BC(\bbhatRl)]$ achieves its minimum,
61.4449, at $\lam$ = 0.0173. In all instances each \emph{VIF}
initially decreases, then increases, but values of $\lam$ at which
the changes occur differ across the five estimators. If we view
$\textup{E}(\bbhatRl)$ = $T(\bbet)$ as a transformation on the
parameter space, Question 4 of Section 4.5 asks whether its
conditioning improves with increasing $\lam.$ To the contrary, the
last column of Table 8 shows that condition numbers increase
explosively with increasing $\lam.$ From Table 3 it is clear that
corresponding condition numbers for $\textup{E}(\bbhatSl)$ =
$T(\bbet)$ are square roots of those listed in Table 8 for
$\bbhatRl.$

Similar entries in Table 9 give the evolution of \emph{VIF}s and
$\ca[\BC(\bbhatSl)]$ for $\bbhatSl.$

A noted departure from Table 8 is that the maximal \emph{VIF} is
$\VM(\bbhatSl)$ = $V\!I\!F(\bhata)$ for all cases, independently of
$\lam.$ Further computations show that the crossing
$\ca[\BC(\bbhatRl)]$ = 99.56217 = $\ca[\BC(\bbhatSl)]$ occurs at
$\lam$ = 0.02750.

\section{Conclusions} \label{S:concl}

Little of the considerable literature on ridge regression is found
to be consistent with the optimization of Hoerl and Kennard (1970)
under equality constraints $\{\bbhat'\bbhat = \csq\},$ and under the
inequality constraints $\{\bbhat'\bbhat\leq\csq\}$ of Balakrishnan
(1963), despite pervasive claims to the contrary.

The problem is traced to (i) a misapplication of LaGrange's
principle; (ii) the false claim that the constrained solutions have
nonsingular distributions, corresponding one--to--one with
$\bbhatLS;$ and (iii) the implied but incorrect assertion that the
ridge parameter $\lam$ corresponds one--to--one with $\csq,$ and
thus the false claim that the solution $\bbhatRl$ of $(\BXpX +
\lam\BIk)\bbet = \BX'\BY$ minimizes the residual sum of squares
among estimators of length $\bbhatRl'\bbhatRl = \cstsq.$ Our Theorem
1 supplies the missing distributions appropriate to constrained
minimization. Generalized ridge regression, seen as solving the
equations $(\BXpX + \Blam)\bbet = \BX'\BY$ with nonnegative ridge
parameters $\Blam = \textup{Diag}(\seq{\lam}{k}),$ is also shown to
be inconsistent with LaGrange minimization.

LaGrange optimization having failed as a rational foundation for
conventional ridge regression, alternatives based on conditioning
are developed in Section 4. Limitations in Type A conditioning, on
which a justification for $\bbhatRl$ rests, prompt the introduction
of \emph{surrogate ridge} solutions, $\bbhatSl,$ to account for
ill--conditioning of $\BX$ on both sides of the \emph{OLS}
equations, $\BXpX\bbet$ = $\BX'\BY.$ Extensive numerical studies, as
reported in Section 5, reexamine the Hospital Manpower Data in a
manner complementary to the conventional analyses undertaken in
Myers (1990). It is demonstrated that none of the conditionings of
$\bbhatRl(\BY),$ $\textup{E}(\bbhatRl) = T(\bbet),$ and
$\textup{V}(\bbhatRl),$ nor the variance inflation factors, as
critical properties of the ridge estimators $\{\bbhatRl;\lam\geq
0\},$ is enhanced monotonically on increasing $\lam.$ In contrast,
for the surrogate solutions $\bbhatSl,$ all (except $T(\bbet)$) of
these are uniformly enhanced as $\lam$ evolves. It is seen that
$\bbhatRl$ is better within a narrow range for small $\lam,$ but its
\emph{VIF}s and condition numbers often become excessive within the
range of $\lam$ often recommended in practice. In short, ridge
regression often exhibits some of the very pathologies it is
intended to redress.

In summary, there is a vast and expanding compendium on the
so--called theory, methodology, and simulation studies surrounding
ridge regression. If indeed constrained optimization is to be
pivotal, then the bulk of these studies will have to be reworked to
take into account the nonstandard distributions of Section 3.3, as
well as constraints for the ridge parameter to be minimizing, as
documented in Sections 3.4 and 5.3. It is remarkable that this field
of applied engineering has thrived for so long, despite critical
false assertions and a dearth of sustaining foundation principles.


\end{document}